\documentclass{conm-p-l}

\usepackage{amsmath}
\usepackage{amsthm}
\usepackage{mathrsfs}
\usepackage{amsfonts}
\usepackage{ amssymb }
\usepackage{amsmath,amscd}
\usepackage{array}
\usepackage{amssymb}
\usepackage[all, cmtip]{xy}
\usepackage{tikz}
\usepackage{tensor}
\usetikzlibrary{arrows}
\usetikzlibrary{cd}

\newtheorem{theorem}{Theorem}[section]

\theoremstyle{definition}
\newtheorem{definition}[theorem]{Definition}
\newtheorem{example}[theorem]{Example}
\newtheorem{stexample}[theorem]{Stack Example}

\theoremstyle{remark}

\numberwithin{equation}{section}

\newcommand{\specialP}{\pmb{\mathbb{P}^1}}

\newcommand{\CC}{\mathbb{C}}
\newcommand{\QQ}{\mathbb{Q}}

\newcommand{\PP}{\mathbb{P}}

\renewcommand{\AA}{\mathbb{A}}

\newcommand{\T}{\mathscr{T}}

\newcommand{\ZZ}{\mathbb{Z}}

\newcommand{\LL}{\mathcal{L}}

\newcommand{\Bun}{\mathcal{B}un}

\newcommand{\OO}{\mathscr{O}}

\newcommand{\M}{\mathcal{M}}
\newcommand{\Mbar}{\mathcal{\overline{M}}}
\newcommand{\MM}{\mathfrak{M}}

\newcommand{\x}{\mathbf{x}}

\renewcommand{\u}{\mathbf{u}}
\renewcommand{\v}{\mathbf{v}}
\renewcommand{\P}{\mathscr{P}}
\newcommand{\Q}{\mathscr{Q}}
\newcommand{\universalP}{\mathfrak{P}}

\newcommand{\g}{\mathfrak{g}}

\newcommand{\bxi}{\pmb{\xi}}

\newcommand{\sslash}{\mathord{/\mkern-6mu/}}

\DeclareMathOperator{\Aut}{Aut}

\DeclareMathOperator{\Hom}{Hom}
\DeclareMathOperator{\Sec}{Sec}

\DeclareMathOperator{\rk}{rk}
\DeclareMathOperator{\vir}{vir}
\DeclareMathOperator{\mov}{mov}
\DeclareMathOperator{\Spec}{Spec}

\begin{document}

\title{Quasimaps and Some Examples of Stacks for Everybody}

\author{Rachel Webb}
\address{Department of Mathematics, University of Michigan, Ann Arbor, Michigan 48109}
\curraddr{Department of Mathematics,
University of California, Berkeley, Berkeley, California, 94720-3840}
\email{webbra@umich.edu}
\thanks{The first author was supported in part by NSF RTG grant 1045119.}

\date{November 30, 2018.}


\begin{abstract}
This note introduces the theory of quasimaps to GIT quotients with intuition and concrete examples, with the goal of explaining a closed formula for the quasimap $I$-function. Along the way, it emphasizes aspects of this story that illustrate general stacky concepts.
\end{abstract}

\maketitle

\section{Introduction}
In modern algebraic geometry, stacks are pervasive. To name just a few appearances, they arise as quotients of varieties by group actions, as the algebraic analog of orbifold curves, and as moduli spaces of objects with nontrivial automorphism groups. Unfortunately the theory is laden with technicalities, making it difficult for newcomers. This note provides concrete examples of stacks that arise in the theory of \textit{quasimaps to GIT quotients.} Quasimaps are themselves a topic in Gromov-Witten theory with a growing list of applications.

Primarily a friendly introduction to quasimaps, this note is also a ``second course'' to Fantechi's brief introduction to stacks \cite{fantechi}, a tour of some examples of stacks ``in the wild.'' 
As an introduction to quasimaps, the goal of this note is to present a closed formula for the quasimap $I$-function. Along the way, we'll highlight examples that illustrate general stacky concepts with the heading {\sc Stack Example}. We will always work over $\CC$.

\section{Motivation for the study of quasimaps}\label{sec:motivation}
The theory of $\epsilon$-stable quasimaps to GIT quotients, as presented in this note, was introduced by Ciocane-Fontanine, Kim, and Maulik in a series of papers: \cite{stable_qmaps}, \cite{wcgis0}, \cite{toric_qmaps}, \cite{bigI}, and \cite{highergenus}. For the sake of novel and concrete exposition, this note omits much of the historical and mathematical context of the theory, which can be found in the original papers and in the survey articles \cite{survey1} and \cite{survey2}.

Gromov-Witten theory begins with the study of maps from a smooth genus-$g$ curve with $n$ marks to a projective target, for example to $\PP^n$. The set of such maps having a fixed degree $d$ forms a moduli space, denoted (in this example) $\M_{g,n}(\PP^n,d)$. Unfortunately, $\M_{g,n}(\PP^n,d)$ is not compact; in order to define Gromov-Witten invariants of $\PP^n$ as integrals over a moduli space of maps, we must replace $\M_{g,n}(\PP^n,d)$ with a compactification. One such compactification is the moduli space of Kontsevich-stable maps $\Mbar_{g,n}(\PP^n,d)$. These stable map moduli spaces are the central objects of Gromov-Witten theory.

\begin{figure}[ht]
\centering
\begin{tikzpicture}[scale=.8]
\draw [dashed] (0,0) rectangle (4,3);
\draw [thick, opacity=.3](1,0)--(1,3);
\draw [thick, opacity=.6](2,0)--(2,3);
\draw [thick] (3,0)--(3,3);
\draw [thick] (2.75,0.65)--(3.5,3.2);
\draw [fill=black] (3,1.5) circle [radius=0.05];

\node at (3,-.35) {$a=0$};
\node at (5.75,2) {\large $\widetilde{ \phi}$};

\draw [ultra thick] (5,1.5)--(6.5,1.5);
\draw [ultra thick] (6.25,1.25)--(6.5,1.5)--(6.25,1.75);

\begin{scope}[shift={(9,1.5)}]
\draw [dashed] (-1.5,-1.5) rectangle (1.5,1.5);
\draw [thick] (-1.5,0)--(1.5,0);
\draw [thick] (0,-1.5)--(0,1.5);
\draw [thick, opacity=.3, domain=-1.5:-.4] plot (\x, {.6/\x});
\draw [thick, opacity=.3, domain=.4:1.5] plot (\x,{.6/\x});
\draw [thick, opacity=.6, domain=-1.5:-.133] plot (\x, {.2/\x});
\draw [thick, opacity=.6, domain=.133:1.5] plot (\x, {.2/\x});
\draw [fill=black] (0,0) circle [radius=0.05];
\end{scope}
\end{tikzpicture}
\caption{The morphism $\widetilde \phi$ maps from $\mathrm{Bl}_{a=y=0}\AA^1\times \PP^1$ on the left to $\PP^2$ on the right, depicted here in the chart of $\PP^2$ where the middle coordinate is nonzero. The varying shades of gray show the fibers of this family of stable maps. The fiber over $a=0$ is a map from a nodal curve.}
\label{fig:stablemaps}
\end{figure}
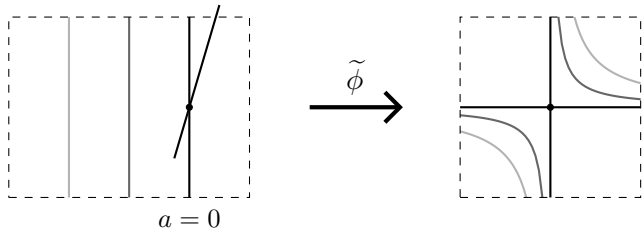

The Kontsevich moduli space compactifies $\M_{g,n}(\PP^n,d)$ by allowing the source curve to be nodal, but there are other ways to compactify. One may see hints of these other ways by looking at possible limits of families of maps in $\M_{g,n}(\PP^n,d)$. Let $C$ be equal to $\PP^1$ with homogeneous coordinates $[x:y]$ and markings at $[1:1]$ and $[2:1]$. Define 
\begin{equation}\label{eq:family}
\phi: \CC^*\times \PP^1 \rightarrow \PP^2 \quad\quad\quad \phi_a(x,y)=[ax^2:xy:y^2] \quad\text{for}\quad a \in \CC^*.
\end{equation}
This is a family in $\M_{0,2}(\PP^2,2)$ with base $\CC^*$. To extend it over the origin to a flat family over $\CC$, we need to define a map $\phi_0$. The natural choice seems to be 
\begin{equation}\label{eq:natural}
\phi_0(s,t)=[0:xy:y^2],
\end{equation}
but this has a basepoint (is undefined) at $y=0$. To recover the limit in $\Mbar_{0,2}(\PP^3,2)$, we resolve the rational map $\phi:\AA^1\times \PP^1\dashrightarrow \PP^2$ (given by \eqref{eq:family} and \eqref{eq:natural}) by blowing up this basepoint, adding an extra rational curve in the fiber over $0$. This produces a morphism $\widetilde \phi$ from the blowup to $\PP^2$, and $\widetilde{\phi}_0$ has a source curve that is two copies of $\PP^1$, glued at a node. The limiting map has degree 1 on each copy. This is depicted in Figure \ref{fig:stablemaps}.

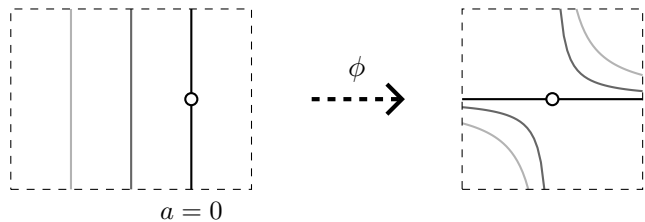
\begin{figure}[ht]
\centering
\begin{tikzpicture}[scale=.8]
\draw [dashed] (0,0) rectangle (4,3);
\draw [thick, opacity=.3](1,0)--(1,3);
\draw [thick, opacity=.6](2,0)--(2,3);
\draw [thick] (3,0)--(3,1.4);
\draw [thick] (3,1.6)--(3,3);
\draw [thick] (3,1.5) circle [radius=.1];

\node at (3,-.35) {$a=0$};
\node at (5.75,2) {\large ${ \phi}$};

\draw [ultra thick, dashed] (5,1.5)--(6.5,1.5);
\draw [ultra thick] (6.25,1.25)--(6.5,1.5)--(6.25,1.75);

\begin{scope}[shift={(9,1.5)}]
\draw [dashed] (-1.5,-1.5) rectangle (1.5,1.5);
\draw [thick] (-1.5,0)--(-.1,0);
\draw [thick] (.1,0)--(1.5,0);
\draw [thick] (0,0) circle [radius=.1];
\draw [thick, opacity=.3, domain=-1.5:-.4] plot (\x, {.6/\x});
\draw [thick, opacity=.3, domain=.4:1.5] plot (\x,{.6/\x});
\draw [thick, opacity=.6, domain=-1.5:-.133] plot (\x, {.2/\x});
\draw [thick, opacity=.6, domain=.133:1.5] plot (\x, {.2/\x});
\end{scope}
\end{tikzpicture}
\caption{The rational map $\phi$ takes $\AA^1\times \PP^1$ on the left to $\PP^2$ on the right. The varying shades of gray show the fibers of this family of quasimaps. In particular, the fiber over $a=0$ is a rational map of degree 1 with a basepoint of length 1 (see Definition \ref{def:stability}).}
\label{fig:quasimaps}
\end{figure}

However, what happens if we compactify $\M_{0,2}(\PP^2,2)$ by allowing basepoints? That is, what if we take $\phi_0(x,y)=[0:xy:y^2]$ as a rational map from $\PP^1$ to $\PP^2$ to be the limit of \eqref{eq:family}? Indeed, this rational map is a \textit{stable quasimap}, depicted in Figure \ref{fig:quasimaps}.

We will see that with the right definitions, stable quasimaps form a moduli space that is just as well-behaved as the Kontsevich moduli spaces. In fact, stability for quasimaps depends on a postive rational parameter $\epsilon$, giving us a whole collection of moduli spaces! This collection has the following advantages:
\begin{enumerate}
\item When $\epsilon > 2$, the quasimap moduli space is equal to the familiar Kontsevich moduli space, and its invariants are Gromov-Witten invariants.
\item When $\epsilon$ is sufficiently small, certain quasimap invariants (the genus-0 invariants) are easier to compute. A generating function for these explicit invariants is called the \textit{quasimap I-function}.
\item One can ``cross the wall,'' relating quasimap invariants for differing values of $\epsilon$, thereby (roughly) expressing Gromov-Witten invariants in terms of the quasimap $I$-function. 
\end{enumerate}
These statements are a heuristic only; for careful statements and their proofs, see \cite{wcgis0} and \cite{highergenus}.
For certain targets, invariants with small $\epsilon$ and $g=0$ were first studied by Givental \cite{giventaleq}, while the strategy in (3) above for computing Gromov-Witten invariants was first employed by Bertram \cite{bertram}.


\section{Quotients in algebraic geometry}

When a reductive algebraic group $G$ acts on an affine variety $X$, we'd like to take the quotient, producing an algebraic object $X/G$. Unfortunately, there may not be a scheme $X/G$ satisfying the universal property of a quotient---if it exists, this scheme is called a \textit{categorical quotient}. Even if a categorical quotient scheme exists, it may not have the expected topology (e.g., be a \textit{geometric quotient}). Therefore, defining a quotient either requires us to modify the original data $X$ and $G$, or to leave the category of schemes. Both strategies arise in quasimap theory. The resulting quotients are respectively called \textit{geometric invariant theory (GIT) quotients} and \textit{stack quotients.}

We briefly summarize the key definitions of GIT for affine $X$ when $G$ acts with no kernel, as found in \cite{king}. Fix a character $\theta$ of $G$. Given a nonnegative integer $n$, a function $f \in \Gamma(X, \OO_X)$ is a \textit{relative invariant of weight} $\theta^n$ if for every $x \in X$ we have $f(g\cdot x) = \theta(g)^nf(x).$
\begin{definition}
A point $x \in X$ is $\theta$\textit{-semistable} if there exists an integer $n\geq 1$ and a relative invariant $f$ of weight $\theta^n$ such that $f(x) \neq 0$. If moreover the dimension of the orbit $G \cdot x$ is equal to the dimension of $G$ and the $G$-action on $\{y \in X \mid f(y) \neq 0\}$ is closed, then $x$ is $\theta$\textit{-stable}.
\end{definition}

The upshot of these definitions is that we get an open locus $X^{ss}_\theta(G)\subset X$ of {semistable points} with respect to $\theta$, and a smaller locus $X^s_\theta(G) \subset X^{ss}_\theta(G)$ of stable points. We define the locus of \textit{unstable points} to be $X^{us}_\theta(G) = X\setminus X^{ss}_\theta(G)$. The sets $X^s_\theta(G)$ and $X^{ss}_\theta(G)$ can be computed with \cite[Prop~2.5]{king}. We include the group $G$ in our notation here because later we will vary it; when there is no risk of misunderstanding, we will omit both the character and the group from the notation. 
In this note, we will assume
\begin{itemize}
\item $X^{ss}_\theta(G)=X^s_\theta(G)$ and this set is nonempty and smooth
\item $G$ acts freely on $X^{ss}_\theta(G)$
\end{itemize}
In this case, the GIT quotient is the smooth variety
\[
X\sslash_\theta G := X^s_\theta(G)/G
\]
where the right hand side is defined to be the categorical quotient (it is a theorem that this exists as a scheme). As a topological space, $X\sslash_\theta G$ is the usual topological quotient $X^s_\theta(G)/G$. Again, when $\theta$ is understood, we will omit it from the notation. In this note, we will also assume that $X\sslash G$ is a projective variety. In general it is projective over the affine quotient $\Spec(\Gamma(X, \OO_X)^G)$.

Whereas the GIT quotient ``forgets'' the unstable locus $X^{us}$, this information is retained in the stack quotient. The objects of the stack quotient $[X/G]$ over a scheme $S$ are diagrams
\[
\begin{tikzcd}
\P \arrow[r] \arrow[d]& X\\
S &
\end{tikzcd}
\]
where $\P$ is a principal $G$-bundle (locally trivial in the \'etale topology) and $\P\rightarrow X$ is a $G$-equivariant map. Morphisms in this category are given by fiber diagrams.

One should picture the GIT quotient $X\sslash G$ as an open subset of the stack quotient $[X/G]$. In a moment we will define quasimaps to $X\sslash G$ to be certain ``stable'' maps to the stack quotient $[X/G]$. Motivation for this definition comes from the previous section and the next example.

\begin{stexample}[GIT vs stack quotient] Let $\CC^*$ act on $\CC^{n+1}$ by
\begin{equation}\label{eq:pnaction}
\lambda \cdot (x_1, \ldots, x_n)=(\lambda x_1, \ldots, \lambda x_n).
\end{equation}
Then if $\theta:\CC^* \rightarrow \CC^*$ is the identity, we have $(\CC^{n+1})^{ss}=\CC^{n+1}\setminus\{0\}$ and $\CC^{n+1}\sslash_\theta \CC^*=\PP^n$. A map from a scheme $S$ to $\PP^n$ is given by a line bundle $\mathcal{L}$ on $S$ and $n+1$ sections of $\mathcal{L}$ that do not simultaneously vanish.

On the other hand, the stack $[\CC^{n+1}/\CC^*]$ ``remembers the origin.'' By definition, a map from a scheme $S$ to $[\CC^{n+1}/\CC^*]$ is a principal $\CC^*$-bundle $\P$ on $S$ and an equivariant map $\P \rightarrow \CC^{n+1}$. As we will see in Example \ref{ex:section}, this equivariant map is equivalent to a section of $\P\times_{\CC^*}\CC^{n+1} = \mathcal{L}^{\oplus n+1} \rightarrow S$ where $\mathcal{L}$ is the line bundle associated to $\P$. Comparing this to the data of a map to $\CC^{n+1}\sslash \CC^*$, we see that the only difference is that now, the sections of $\mathcal{L}$ are allowed to vanish simultaneously---i.e., we allow the map to ``hit the origin.''

\end{stexample}

We close this section with an important example in which a desirable quotient does exist as a scheme, namely certain \textit{mixing spaces}. Let $\P\rightarrow S$ be a principal $G$-bundle with a left $G$-action and let $X$ be an affine $G$-variety. Then $\P \times X$ carries a natural $G$-action defined on closed points by
\[
g\cdot(p,x) = (gp, gx).
\]
We define the \textit{mixing space} to be the categorical quotient
\[
\P \times_G X := (\P \times X) /G.
\]
 One can show that the categorical quotient exists as a scheme, and in fact it is a geometric quotient (see for example the discussion in \cite[Sec~3]{brion}). It has a natural map to $\P/G = S$, with closed fibers isomorphic to $X$. In particular, when $X$ is a vector space and $G$ acts linearly on $X$, the mixing space $\P\times_G V$ is the total space of a vector bundle on $S$.

In fact, there is a bijection between principal $GL_r$-bundles and rank-$r$ vector bundles on $S$ given by sending $\P$ to $\P\times_G \CC^r$, where $\CC^r$ is a $GL_r$-module via left multiplication. We call $\P\times_G \CC^r$ the \textit{vector bundle associated to }$\P$, and $\P$ the \textit{underlying principal bundle of} $\P\times_G \CC^r$.

\section{Quasimaps}
This note will discuss the theory of quasimaps to $V\sslash G$ when $V$ is a vector space. (The original papers present this theory in more generality, allowing $V$ to be an affine variety with l.c.i. singularities contained in $V^{us}$.) Let $V$ be a vector space and $G$ a reductive algebraic group that acts linearly from the left on $V$. We've seen that a quasimap to $\PP^2$ is, loosely speaking, a rational map to $\PP^2$, and that such a map can be realized as a map to the stack quotient $[\CC^3/\CC^*]$. This motivates the following definition, which essentially says that a prestable quasimap to $V\sslash G$ is a map to the stack quotient $[V/G]$ that recovers a rational map.

\begin{definition}
A prestable \textit{quasimap} to $V\sslash G$ is data $(C, x_1, \ldots, x_n, \P, \tilde u)$ where
\begin{itemize}
\item $C$ is a curve with at worst nodal singularities and marked points $x_1, \ldots, x_n$
\item $\varrho: \P\rightarrow C$ is a principal $G$-bundle with a left $G$-action 
\item $\tilde u: \P \rightarrow V$ is a $G$-equivariant map with $\varrho(\tilde u^{-1}(V^{us}))$ a finite set that is disjoint from nodes and markings.
\end{itemize}
The \textit{genus} of the quasimap is the genus of $C$, and the \textit{degree} of the quasimap is the degree of $\P$, i.e., a group homomorphism $\beta \in Hom(\chi(G),\ZZ)$ defined by 
\[
\beta(\omega)=\deg_C(\P\times_G \CC_\omega)
\]
where $\omega \in \chi(G)$ is a character of $G$ and $\CC_\omega$ is the corresponding 1-dimensional representation
\[
g\cdot z = \omega(g)z \quad \quad \text{for}\;g \in G\;\text{and}\;z \in \CC_\omega.
\]
\end{definition}

To write down examples of quasimaps, it is helpful to replace $\tilde u$ with an associated section of $\P\times_G V \rightarrow C$, which we will typically denote $u$. In fact this is equivalent data, as explained in the example below.

\begin{stexample}[Morphisms of stacks]\label{ex:section} If $\varrho: \P \rightarrow C$ is a principal $G$-bundle, then equivariant maps $\P \rightarrow V$ are in bijection with sections of the associated vector bundle $\P\times_G V \rightarrow C$. I'll sketch this bijection.

From the universal property of fiber products, we have a natural bijection between $\Hom(\P, V)$ and $\Sec(\P,\P\times V)$ where $\Sec$ denotes the space of sections. Letting $G$ act on morphisms by conjugation, this is a $G$-equivariant bijection, so we have 
\[
\Hom_G(\P,V) \cong \Sec_G(\P, \P\times V).
\]
where the subscript $G$ indicates the $G$-invariant part (in particular, $\Hom_G(\P,V)$ is the set of $G$-equivariant maps). On the other hand, I'll sketch a bijection
\[
\Sec(C, \P\times_G V) \cong \Sec_G(\P, \P\times V)
\]
by interpreting $\P\times_G V$ as the stack quotient $[(\P\times V)/G]$. A map from $C$ to $\P\times_G V$ is by definition a principal $G$-bundle $\Q$ on $C$ and an equivariant map to $\P\times V$; to say that it is a section means that composition with projection to $C$ is the identity:
\[
\begin{tikzcd}
\Q \arrow[r] \arrow[d] & \P\times V \arrow[r] & \P \arrow[d]\\
C \arrow[rr, equal] && C
\end{tikzcd}
\]
Then the induced map $\Q \rightarrow \P$ is a morphism of principal bundles, hence an isomorphism. After identifying $\P$ with $\Q$ this way, the map $\Q \rightarrow \P\times V$ becomes a $G$-equivariant section $\P \rightarrow \P\times V$. The reader is invited to find the inverse to this correspondence.
\end{stexample}

\begin{example}\label{ex:qmaps}
For two GIT quotients $V\sslash G$, we'll write down all quasimaps from $\PP^1$ to $V\sslash G$ as vectors of homogeneous polynomials. This example is generalized in Section \ref{sec:mainthm}.
\begin{enumerate}
\item If $V=\CC^{n+1}$ and $G = \CC^*$ with action \eqref{eq:pnaction}, and if $\theta$ is the identity character of $\CC^*$, then $V^{ss}$ is $V \setminus \{0\}$ and $V\sslash_{\theta} G$ is $\PP^n$. A quasimap to $\PP^n$ of degree $d \in \Hom(\ZZ, \ZZ)=\ZZ$ with source curve $\PP^1$ is
\begin{itemize}
\item A principal $\CC^*$-bundle $\P$ such that $deg_{\PP^1}(\P\times_{\CC^*}\CC_{\theta})=d$. 
Therefore $\P$ is the underlying principal bundle of $\OO_{\PP^1}(d).$
\item A section of $\P\times_{\CC^*}\CC^{n+1}.$ Here the action of $\CC^*$ on $\CC^{n+1}$ is given by \eqref{eq:pnaction}, so that this vector bundle is $n+1$ copies of the associated bundle to $\P$, i.e., $\P\times_{\CC^*}\CC^{n+1}$ equals $\OO(d)^{\oplus n+1}$.
\end{itemize}
Therefore, a prestable quasimap to $\PP^n$ of degree $d \in \ZZ$ is given by a vector
\[
(\;p_1(x,y), \;p_2(x,y), \;\ldots,\; p_{n+1}(x,y)\;)
\]
of homogeneous polynomials of degree $d$ that are not all zero. In particular \eqref{eq:natural} is a prestable quasimap of degree 2.
\item If $V=M_{k\times n}$ is $k\times n$ matrices over $\CC$ and $G=GL_k$ acts on $V$ by left multiplication and $\theta$ is the determinant character, then $V^{ss}$ is full-rank matrices and $V\sslash_{\theta} G$ is the Grassmannian $Gr(k,n)$. A quasimap to $Gr(k,n)$ of degree $d \in \Hom(\ZZ,\ZZ)=\ZZ$ with source curve $\PP^1$ is
\begin{itemize}
\item A principal $GL_k$-bundle $\P$ on $\PP^1$ such that $\deg_{\PP^1}(\P\times_{GL_k}\CC_\theta) = d$. 
Let $E$ denote the vector bundle associated to $\P$.
From Grothendieck's classification of principal bundles \cite{gro57}, we have $E = \oplus_{i=1}^k \OO(d_i)$ for some $d_i$ with $\sum_{i=1}^k d_i=d$. 
\item A section of $\P\times_{GL_k}M_{k\times n}$, i.e., of $E^{\oplus n}$.
\end{itemize}
So a prestable quasimap to $Gr(k,n)$ of degree $d$ is given by a matrix of polynomials
\begin{equation}\label{eq:matrix}
[p_{ij}(x,y)]_{1\leq i \leq k, \; 1\leq j \leq n}
\end{equation}
where $p_{ij}(x,y)$ is homogeneous of degree $d_i$. Because $M_{k\times n}^{us}$ is matrices of low rank, to define a prestable quasimap the matrix \eqref{eq:matrix} must have low rank on a finite set.
\end{enumerate}
\end{example}

\section{Stability and moduli spaces of quasimaps}
In analogy with Gromov-Witten invariants, we want to define quasimap invariants of a target $V\sslash G$ to be integrals on certain moduli spaces. This cues the entrance of our second example of a stack: the moduli space of quasimaps. We will see that this is a ``master'' moduli space containing $\M_{g,n}(V\sslash G, \beta)$ and many compactifications of it. Then we will define some of these compactifications (substacks of the master space) to be the $\epsilon$-stable moduli spaces of quasimaps.

\begin{stexample}[Moduli stack of prestable quasimaps]
Let $\MM_{g,n}(V\sslash G, \beta)$ be the moduli stack of prestable quasimaps of genus $g$, degree $\beta$, and $n$ marks. Objects in this category over a base $S$ are families of prestable quasimaps on $S$; i.e., they 
are triples $(\mathcal{C}, \P, \tilde u)$
where $\mathcal{C}\rightarrow S$ is a flat family of genus-$g$ nodal curves (not necessarily stable) on $S$, $\P\rightarrow \mathcal{C}$ is a principal $G$-bundle, $\tilde u:\P \rightarrow V$ is $G$-equivariant, and geometric fibers over $S$ are prestable quasimaps of degree $\beta$. An isomorphism between objects in $\MM_{g,n}(V\sslash G, \beta)(S)$ is a commuting diagram
\[
\begin{tikzcd}
\P' \arrow[r, "\sim"]\arrow[rr, bend left, "\tilde u'"] \arrow[d] &\P\arrow[d] \arrow[r,"\tilde u"] & V\\
\mathcal{C}'\arrow[r,"\sim"]&\mathcal{C}
\end{tikzcd}
\]
where $\mathcal{C}' \xrightarrow{\sim} \mathcal{C}$ commutes with the maps to $S$ and the square is fibered (such a diagram is an isomorphisms of quasimap families). The stack $\MM_{g,n}(V\sslash G, \beta)$ is not Deligne-Mumford, as some prestable quasimaps have non-finite automorphism groups. We see the familiar offenders from stable map theory: for example, a degree-0 map sending $\PP^1$ to a point in $\PP^2$ is invariant under the entire automorphism group of $\PP^1$. However, there are new examples as well. For instance, 
$[x^2: 3x^2: x^2]$
defines a prestable quasimap of degree 2 from $\PP^1$ to $\PP^2$ which is invariant under $[x:y]\mapsto[x:ty]$.
\end{stexample}

\begin{example}\label{ex:qmap_isos} Let's find isomorphisms between the quasimaps in Example \ref{ex:qmaps}.
\begin{enumerate}
\item Let $V, G$, and $\theta$ be as in Example \ref{ex:qmaps} part 1. From that example, a prestable quasimap from $\PP^1$ to $\PP^n$ of degree $d$ is given by a vector of homogeneous degree-$d$ polynomials $(p_i(x,y))_{i=1}^{n+1}$. Then a quasimap isomorphism is an element $\alpha \in \Aut(\OO(d))=\CC^*$, and it sends $(p_i(x,y))$ to $(\alpha p_i(x,y))$.
\item Let $V, G$, and $\theta$ be as in Example \ref{ex:qmaps} part 2. From that example, a prestable quasimap from $\PP^1$ to $Gr(k,n)$ of degree $d$ is $n$ sections of a vector bundle $\oplus_{i=1}^k \OO(d_i)$ of degree $d$, which may be denoted by a $k\times n$ matrix $(p_{ij}(x,y))$ of polynomials where $p_{ij}(x,y)$ is homogeneous of degree $d_i$. Then a quasimap isomorphism is an element $A$ of $\Aut(E) = \Hom(\oplus \OO(d_i), \oplus \OO(d_i))^\times$ which we may identify with a $k\times k$ matrix $(a_{\ell i}(x,y))$ of polynomials where $a_{\ell i}(x,y)$ is homogeneous of degree $d_\ell-d_i$. Such an isomorphism acts on a the quasimap $(p_{ij}(x,y))$ by matrix multiplication. (Notice that if $d_1\geq d_2\geq \ldots \geq d_k$, then $(a_{\ell i}(x,y))$ will be block upper triangular.)
\end{enumerate}
\end{example}

For many reasons, the stack $\MM_{g,n}(V\sslash G,\beta)$ is not the right one for defining invariants. Instead, we impose a \textit{stability condition} that cuts out a proper separated Deligne-Mumford substack of $\MM_{g,n}(V\sslash G, \beta)$ which contains $\M_{g,n}(V\sslash G, \beta)$ as an open subset. 

\begin{definition}\label{def:stability}
Let $(C, x_1, \ldots, x_n, \P, \tilde u)$ be a prestable quasimap to $V\sslash_\theta G$ and let $L = \P\times_G \CC_\theta$. This quasimap is $\epsilon$-stable if
\begin{enumerate}
\item On every component $C'$ of $C$ we have
\[
2g_{C'}-2+n_{C'}+\epsilon \deg(L|_{C'}) > 0
\]
where $g_{C'}$ is the genus of $C'$ and $n_{C'}$ is the number of marked points and nodes on $C'$, and
\item For every $x \in C$ we have
\[
\ell(x) \leq 1/\epsilon
\]
where $\ell(x)$, called the \textit{length} of $x$, is the order of contact of $u(C)$ with $\P\times_G V^{us}$, and $u$ is the section of $\P\times_G V\rightarrow C$ determined by $\tilde u$. See \cite[Def~7.1.1]{stable_qmaps} for more on the definition of length. For $x \in C$, the length $\ell(x)$ is nonzero if and only if $u(x)$ is in $\P\times_G V^{us}$, and in this case we say $x$ is a \textit{basepoint} of the quasimap.
\end{enumerate}
\end{definition}

The first condition (1) is sometimes stated as ``$\omega_C(\sum x_i)\otimes L^\epsilon$ is ample.'' A family of prestable quasimaps is $\epsilon$-stable if every geometric fiber is $\epsilon$-stable. The moduli space of $\epsilon$-stable quasimaps is denoted $\Mbar_{g,n}^\epsilon(V\sslash G, \beta)$. An $\epsilon$-stable quasimap only has finitely many automorphisms. Hence the moduli spaces $\Mbar_{g,n}^\epsilon(V\sslash G, \beta)$ are Deligne-Mumford, not just Artin stacks. They also have other good geometric properties, as stated in the following theorem.

\begin{theorem}\cite[Thm~7.1.6]{stable_qmaps}\label{thm:dmstack}
The moduli space $\Mbar_{g,n}^\epsilon(V\sslash G, \beta)$ is a proper separated Deligne-Mumford stack of finite type.
\end{theorem}

As explained in the introduction, the benefit of these spaces, especially in genus 0, is that when $\epsilon$ is greater than 2, the space $\Mbar_{0,n}^\epsilon(V\sslash G,\beta)$ is the familiar space $\Mbar_{0,n}(V\sslash G, \beta)$. On the other hand, when $\epsilon$ is sufficiently small, the moduli spaces $\Mbar_{0,n}^\epsilon(V\sslash G,\beta)$ do not depend on $\epsilon$; this is called the $0+$-stable quasimap moduli space. The space $\Mbar_{0,n}^{0+}(V\sslash G, \beta)$ is more ``computable'' (see Example \ref{ex:represent}). Hence the wall-crossing theorem of \cite{wcgis0}, which translates between invariants of $\Mbar_{0,n}^\epsilon(V\sslash G, \beta)$ for differing values of $\epsilon$, gives a way to relate the Gromov-Witten invariants of $V\sslash G$ to invariants that are more computable.

\begin{example}
In Section \ref{sec:motivation} we described two possible limits of the family \eqref{eq:family} of quasimaps to $\PP^2$. One limit was the stable map $\widetilde{\phi}$. Because $\widetilde{\phi}$ has no basepoints, it satisfies condition (2) in Definition \ref{def:stability} for each $\epsilon$. However, on the component with no marks---call it $C'$---we have
\[
2g_{C'}-2+n_{C'}+\epsilon \deg(L|_{C'})=-1+\epsilon,
\]
which is possible only for $\epsilon > 1$. So this map is $\epsilon$ stable for $\epsilon > 1$.

The second limit was the rational map $\phi_0(x,y)=[0:xy:y^2]$. This quasimap satisfies condition (1) of Definition \ref{def:stability} for every $\epsilon >0$. However, we have $\ell([1:0])=1$, so that this map is $\epsilon$-stable only when $\epsilon \leq 1$.

Hence, the family \eqref{eq:family} in $\Mbar_{g,n}^\epsilon(V\sslash G, \beta)$ has the stable map $\widetilde{\phi}$ for a limit when $\epsilon > 1$, and the rational map $\phi_0$ for a limit when $\epsilon \leq 1$.
\end{example}

A variant of $\Mbar^\epsilon_{g,n}(V\sslash G, \beta)$ are the $\epsilon$-stable \textit{quasimap graph spaces}. Graph spaces are used in many contexts to prove wall crossing or mirror theorems. In particular, we will use the graph moduli space with $g=n=0$ and $\epsilon=0^+$ in Section \ref{sec:Ifunc} to define the $I$-function, so we'll define the quasimap graph space with these parameters now. 

\begin{definition}\label{def:QG}Fix a ``reference copy'' of $\PP^1$ and denote it $\specialP$. The quasimap graph space $QG(V\sslash G, \beta)$ is the stack whose objects over a scheme $S$ are prestable quasimaps $(\mathcal{C}, \P, \tilde u)$, with an additional datum
\[
\phi: \mathcal{C}\rightarrow \specialP
\]
which restricts to an isomorphism on every geometric fiber of $\mathcal{C}\rightarrow S$.
An isomorphism between two objects $(\mathcal{C}, \P, \tilde u, \phi)$ and $(\mathcal{C}',\P',\tilde u', \phi')$ in $QG(V\sslash G, \beta)(S)$ is a commuting diagram
\begin{equation}\label{fig:qgmorph}
\begin{tikzcd}
&\P' \arrow[r, "\sim"]\arrow[rr, bend left, "\tilde u'"] \arrow[d] &\P\arrow[d] \arrow[r,"\tilde u"] & V\\
\specialP & \arrow[l,"\phi'"'] \mathcal{C}'\arrow[r,"\sim"]&\mathcal{C} \arrow[ll, bend left, "\phi"]
\end{tikzcd}
\end{equation}
such that $\mathcal{C}' \xrightarrow{\sim} \mathcal{C}$ commutes with the maps to $S$ and the square is fibered.
\end{definition}
When $S$ is $\mathrm{Spec}(\CC)$, The additional datum $\phi$ realizes $\mathcal{C}$ as $\specialP$. So closed points of $QG(V\sslash G, \beta)$ are in bijection with quasimaps from $\specialP$ to $V\sslash G$. In addition, notice that $\phi$ prevents the existence of automorphisms, so we do not need marked points to achieve stability. According to \cite[Thm~7.2.2]{stable_qmaps}, the space $QG(V\sslash G,d)$ is a proper Deligne-Mumford stack of finite type. In the next example, we identify $QG(\PP^n, d)$ as a smooth scheme.

\begin{stexample}[A moduli stack represented by a scheme]\label{ex:represent}
We have seen (Example \ref{ex:qmaps}) that a quasimap from $\PP^1$ to $\PP^n$ of degree $d$ is a nonzero element of $\Gamma( \PP^1, \OO(d))^{\oplus n+1}$. However, two such sections define isomorphic quasimaps exactly when they differ by a complex scalar (Example \ref{ex:qmap_isos}). Hence, we naively expect 
\[QG(\PP^n,d) = \left(\;\Gamma(\PP^1, \OO(d))^{\oplus n+1}\;\setminus \;\{0\}\;\right)/\;\CC^*\cong \PP^N,\] where $N=dn+d+n$.
Indeed, this space carries a tautological family of quasimaps as follows. On $\PP^N$ we have the trivial family of curves $\PP^N \times \PP^1$, and on this family the vector bundle $V = \OO_{\PP^N}(1)^{\oplus n+1}\otimes \OO_{\PP^1}(d)$. An element of $\PP^N \times \PP^1$ may be written $(\sigma, \x)$ where $\sigma \in \Gamma( \PP^1, \OO(d))^{\oplus n+1}$ is a vector of $n+1$ degree-$d$ homogeneous polynomials in two variables and $\x = (x,y)$. The tautological quasimap is given by the section of $V$ sending $(\sigma, \x)$ to $\sigma(\x)$.

The tautological family on $\PP^N$ defines a map $F:\PP^N \rightarrow QG(\PP^n, d)$ which is a bijection on closed points by Examples \ref{ex:qmaps} and \ref{ex:qmap_isos}. It can be shown that $QG(\PP^n, d)$ is a smooth algebraic space (see Examples \ref{ex:obs_vanish} and \ref{ex:algspace}). Since $\PP^N$ is smooth as well, it follows that $F$ must be an isomorphism (see \cite[Lem~3.6.2]{dissertation}).
\end{stexample}

\section{Virtual cycles via perfect obstruction theories}
We want to define invariants of $X\sslash G$ which are integrals over the moduli spaces\\ $\Mbar^\epsilon_{g,n}(X\sslash G, \beta)$. But these spaces may not be smooth or even equidimensional, so the usual definition of integration doesn't make sense. We need to choose some class to be the (virtual) fundamental class with which to define integration. One way to get such a class is via a perfect obstruction theory and the intrinsic normal cone construction of Behrend-Fantechi \cite{bf}. We will summarize their construction; for a more thorough but still elementary introduction see \cite{bfm}.

Let $X$ be a Deligne-Mumford stack with a map to an Artin stack $Y$. We'll assume $X$ and $Y$ are separated and of locally finite type, and that the map $X \rightarrow Y$ is of \textit{relative Deligne-Mumford type}, a technical condition that is always satisfied in our examples. In this case, the \textit{relative cotangent complex} of $X$ over $Y$ is a certain object $\LL_{X/Y}^{\bullet}$ in the derived category $D^{\leq 0}(X)$ of quasi-coherent sheaves. One property of $\LL^{\bullet}_{X/Y}$ is that $H^0(\LL_{X/Y}^\bullet)$ is isomorphic to the sheaf of Kahler differentials $\Omega^1_{X/Y}$.

A \textit{relative perfect obstruction theory} on $X$ is a perfect complex $E^\bullet \in D^{[-1,0]}(X)$ and a morphism $\phi: E^\bullet \rightarrow \LL_{X/Y}^{\bullet}$
such that $h^0(\phi)$ is an isomorphism and $h^{-1}(\phi)$ is surjective. From this data, Behrend-Fantechi construct
\begin{itemize}
\item A \textit{relative intrinsic normal cone} $\mathfrak{C}_{X/Y}$ that is an Artin stack over $X$
\item A cone stack $\mathfrak{E}\rightarrow X$ that contains the intrinsic normal cone via a closed embedding $\mathfrak{C}_{X/Y} \hookrightarrow \mathfrak{E}$
\end{itemize}
By the second bullet, the intrinsic normal cone is a closed substack of $\mathfrak{E}$; intuitively, the virtual class induced by $E^\bullet$ is the intersection of the intrinsic normal cone with the 0-section of $\mathfrak{E}$. Symbolically,
\[
[X]^{\vir} := 0^!_{\mathfrak{E}}[\mathfrak{C}_{X/Y}]\in A_{\dim Y+\rk E^0-\rk E^{-1}}(X).
\]
The symbol $0^!_{\mathfrak{E}}$ is defined in \cite{bf}.
One incidental use of perfect obstruction theories is to detect when a stack is smooth, using the following example (see \cite[Prop~7.3]{bf}).

\begin{example}\label{ex:obs_vanish}
Let $X$ and $Y$ be as above and let $\phi: E^\bullet \rightarrow \LL_{X/Y}^{\bullet}$ be a relative perfect obstruction theory. If $E^{-1}=0$ then we say \textit{the obstructions vanish.} If this is the case and moreover $E^0$ is locally free, then $h^0(E^\bullet)=E^0$ is locally free, but $h^0(\phi)$ is an isomorphism, so $h^0(\LL_{X/Y}^{\bullet})=\Omega^1_{X/Y}$ is locally free and $X \rightarrow Y$ is smooth. In this case, the construction of Behrend-Fantechi returns the usual fundamental class of $X$; i.e., $[X]^{\vir}=[X]$.
\end{example}

There is a natural relative perfect obstruction theory on $\Mbar_{g,n}^\epsilon(V\sslash G, \beta)$. This stack has a forgetful map to the moduli stack of principal $G$-bundles on genus $g$ curves
\[
\mu: \Mbar_{g,n}^\epsilon(V\sslash G, \beta) \rightarrow \Bun_G
\]
given by forgetting the section from the quasimap data. 
\begin{theorem}\cite[Thm~7.1.6]{stable_qmaps}
On $\Mbar_{g,n}^\epsilon(V\sslash G, \beta)$ 
let $\pi: \mathcal{C} \rightarrow \Mbar_{g,n}^\epsilon(V\sslash G, \beta)$ 
be the universal curve, let 
$\universalP$ denote the universal principal bundle, and let 
$\varrho: \universalP\times_G V \rightarrow \mathcal{C}$ be the associated vector bundle with 
$T_\varrho$ the relative tangent bundle of this map. Then the complex
\[
\big( R^\bullet \pi_*(u^*T_\varrho) \big)^\vee
\]
is a $\mu$-relative perfect obstruction theory for $\Mbar_{g,n}^\epsilon(V\sslash G, \beta)$.
\end{theorem}

\noindent
According to \cite{stable_qmaps}, the analogous result holds for the quasimap graph space $QG(V\sslash G, \beta)$. 

\begin{stexample}[Perfect obstruction theory]
We compute the relative perfect obstruction theory on $QG(\PP^n, d)$ and show that it recovers the usual fundamental class of $\PP^N$. From Stack Example \ref{ex:represent}, the relative perfect obstruction theory is the complex
\[
\big( R^\bullet \pi_*(u^*T_\varrho) \big)^\vee = \big(R^\bullet \pi_*
( \OO_{\PP^1}(d)\otimes \OO_{\PP^N}(1)^{\oplus n+1})
\big)^\vee
\]
where $\pi: \PP^N \times \PP^1 \rightarrow \PP^N$ is projection to the first factor. By the projection formula, this is the dual of the complex
\[
R^0\pi_*(\OO_{\PP^1}(d))\otimes \OO_{\PP^N}(1)^{\oplus n+1} \rightarrow R^1\pi_*(\OO_{\PP^1}(d))\otimes \OO_{\PP^N}(1)^{\oplus n+1}.
\]
The second term vanishes, so the $\mu$-relative perfect obstruction theory is the vector bundle $\OO_{\PP^N}^{\oplus d+1}\otimes \OO_{\PP^N}(1)^{\oplus n+1}=\OO_{\PP^N}(1)^{\oplus N+1}$ in degree 0. By Example \ref{ex:obs_vanish} the virtual class defined by this theory is the usual fundamental class on $\PP^N$. 
\end{stexample}

\section{Quasimap $I$-functions}\label{sec:Ifunc}
\sloppy For any positive rational $\epsilon$, we have a proper separated Deligne-Mumford stack $\Mbar_{0,n}^\epsilon(V\sslash G, \beta)$ with a perfect obstruction theory. This is exactly what we need to define invariants of $V\sslash G$ as integrals on these spaces. As discussed at the end of Section \ref{sec:motivation}, when $\epsilon$ is at least 2, these invariants are equal to Gromov-Witten invariants. The mirror theorem of \cite{wcgis0} relates these invariants to localization residues on the \textit{graph moduli space} $QG(V\sslash G, \beta)$ for the $\epsilon=0^+$ stability parameter (see Definition \ref{def:QG}). The residues on $QG(V\sslash G, \beta)$ are indexed by the \textit{quasimap $I$-function}. In contrast to other Gromov-Witten generating functions, the advantage of the $I$-function is that it can be written down with a closed formula. The goal of this section is to define the quasimap $I$-function as a formal power series indexing certain localization residues on $QG(V\sslash G, \beta)$.


We already know from \cite[Thm~7.2.2]{stable_qmaps} that $QG(V\sslash G, \beta)$ is a Deligne-Mumford stack, but in fact more is true: it is an \textit{algebraic space}.

\begin{stexample}[$QG(V\sslash G, \beta)$ is an algebraic space]\label{ex:algspace} By \cite[Thm~2.2.5]{conrad} it suffices to show that a quasimap in this space has a trivial automorphism group. Essentially, this is because quasimaps are required to map into $V^s$ generically, and $G$ acts on $V^s$ with trivial stabilizers by our assumptions in Section \ref{sec:motivation}. Let $\sigma: \PP^1 \rightarrow \P\times_G V$ be a quasimap, and suppose $\phi$ is an automorphism of $\P$ commuting with $\sigma$. In a local chart $U$ where $\P$ is trivial, let $\sigma_U:U\rightarrow V$ and $\phi_U:U \rightarrow G$ be coordinate representations of $\sigma$ and $\phi$. Then
\[
\phi_U(u)\sigma_U(u)=\sigma_U(u).
\]
For every $u\in U$ that is not a basepoint, we know $\sigma_U(u)\in V^s$, so this equation implies $\phi_U(u)=1$ in $G$. Hence $\phi_U$ is the identity on a dense subset of $U$, hence on all of $U$.
\end{stexample}

To define the coefficients of the $I$-function as localization residues, we need a $\CC^*$-action on $QG(V\sslash G, \beta)$. This space carries a $\CC^*$ action as follows. Let $x_0, x_1$ be homogeneous coordinates on $\specialP$ and let $\CC^*$ act on $\specialP$ by
\[
\lambda \cdot [x_0: x_1] = [\lambda x_0: x_1], \quad \quad\quad\quad \lambda \in \CC^*.
\]
This induces an action on $QG(V\sslash G, \beta)$ given by
\[
\lambda(\mathcal{C},\P,\tilde u, \phi) = (\mathcal{C},\P,\tilde u, \lambda \circ \phi).
\]
The fixed locus of $\CC^*$ is the closed subspace of $QG(V\sslash G, \beta)$ whose objects over a scheme $S$ are families fixed by the action.

\begin{stexample}[Fixed locus of a group action]
What are the basepoints of a fixed graph quasimap? If a graph quasimap $(\mathcal{C},\P,\tilde u, \phi)$ over $\mathrm{Spec}(\CC)$ is $\CC^*$-fixed, then for every $\lambda \in \CC^*$ we have a diagram \eqref{fig:qgmorph} with $(\mathcal{C}', \P', \tilde u', \phi')=(\mathcal{C},\P,\tilde u, \lambda \circ \phi)$ and $\phi$ an isomorphism. 
Then the map $\mathcal{C}'\rightarrow \mathcal{C}$ in \eqref{fig:qgmorph} must be $\phi^{-1}\circ \lambda \circ \phi$. But $\phi^{-1}\circ \lambda \circ \phi$ must fix basepoints of $\tilde u$, for every $\lambda$. This means that basepoints of $\tilde u$ have to be $\phi^{-1}([0:1])$ or $\phi^{-1}([1:0])$. 

We can use this information to identify components of the fixed locus: the lengths of these basepoints $\ell([0:1])$ and $\ell([1:0])$ are constant in families, so specifying these integers specifies a component of the fixed locus.
\end{stexample}

As explained in the previous example, the fixed locus of $\CC^*$ acting on $QG(V\sslash G, \beta)$ has components determined by the length of the basepoints at $[0:1]$ and $[1:0]$. Let $F_\beta$ be the component where all the basepoints are at $[0:1]\in \specialP$. This component has a natural map $ev_\bullet$ to $V\sslash G$.

\begin{stexample}[Morphism of stacks] Define a map $ev_\bullet: F_\beta \rightarrow V\sslash G$ as follows. Let $(\mathcal{C},\P,\tilde u, \phi)$ be an object of $F_\beta$ lying over $S$. Recall from Example \ref{ex:section} that $\tilde u$ defines a section $u$ of $\P\times_G V \rightarrow \mathcal{C}$; since $[1:0]$ is not a basepoint, $u([1:0])$ is in $\P\times_G V^s$. So define $ev_\bullet$ to send $(\mathcal{C},\P,\tilde u, \phi)$ to the morphism $S\rightarrow V\sslash G$ apparent in the following diagram:
\[
\begin{tikzcd}
&\P\times_G V \arrow[d] & \P\times_G V^s \arrow[l, hook'] \arrow[r] & V^s/G = V\sslash G\\
S\times\{[1:0]\} \arrow[r] & \mathcal{C}\arrow[u,bend left, "u"] &
\end{tikzcd}
\]
\end{stexample}

Finally we can define the $I$-function as follows. The $\mu$-relative perfect obstruction theory on $QG(V\sslash G, \beta)$ 
defines an \textit{absolute perfect obstruction theory} which is $\CC^*$-equivariant. 
The moving part of this equivariant theory may be used to define the euler class of the normal bundle to the fixed locus, denoted
$e_{\CC^*}(N^{\vir}_{F_\beta})$. The fixed part, when restricted to 
$F_\beta$, 
is a perfect obstruction theory for 
$F_\beta$, and hence defines a virtual class 
$[F_\beta]^{\vir}$ (see \cite[Sec~3]{ckl}). With these definitions, we may write the $I$-function 
of $V\sslash G$ as a formal sum: we use formal symbols 
$q^\beta$ to index coefficients that are localization residues on the stacks $F_\beta$. Precisely,
\begin{equation}\label{eq:Ifunc}
I^{V\sslash G}(z) = 1+ \sum_{\beta\neq 0}q^\beta I^{V\sslash G}_\beta(z) \quad \quad \text{where} \quad \quad I^{V\sslash G}_\beta(z) = (ev_{\bullet})_*\left(\frac{[F_\beta]^{vir}}{e_{\CC^*}(N^{vir}_{F_{\beta}})}\right).
\end{equation}

The promise of $0^+$-stable quasimaps was that their invariants would be computable. Indeed, the moduli spaces $F_\beta$ in the formula \eqref{eq:Ifunc} are represented by smooth schemes, as stated in the following theorem:
\begin{theorem}
The moduli space $F_\beta$ is a disjoint union of flag bundles on smooth subvarieties of $V\sslash G$.
\end{theorem}
In \cite{webb} this identification is explicit. The strategy of the proof is the same as that used in Example \ref{ex:represent}: show that $F_\beta$ is a smooth algebraic space, and then find a tautological family on the desired smooth variety (a closed subscheme of a flag bundle on $V\sslash G$) whose geometric fibers are in bijection with objects parametrized by (a component of) $F_\beta$. This theorem can be used to derive an explicit formula for the quasimap $I$-function, a formula we'll write down in the next section.

\begin{example}\label{ex:Fbeta}
Let us describe $F_d$ when the target is the Grassmannian $Gr(k,n)$. Recall from Example \ref{ex:qmaps} that a degree-$d$ graph quasimap to $Gr(k,n)$ is $n$ sections of a vector bundle $\oplus_{\specialP}(\OO(d_i))$ with $\sum d_i=d$, which may be represented by a $k\times n$ matrix of homogeneous polynomials. For convenience assume $d_1 \geq d_2 \geq \ldots \geq d_k$. Clearly, full-rank matrices $(c_{ij}x^i)$ for complex numbers $c_{ij}$ define elements of $F_\beta$. Referring to Example \ref{ex:qmap_isos} we see that invertible matrices of the form $(a_{ij}x^{d_i-d_j})$ define isomorphisms of this collection of elements of $F_\beta$, where $a_{ij}$ are complex numbers with $a_{ij}=0$ if $d_j>d_i$. From this, one may naively guess that there is a component of $F_\beta$ given by
\[
F_{d_1,\ldots, d_k} = \big(M_{k\times n}\setminus \Delta \big)/ U
\]
where $\Delta \subset M_{k\times n}$ is matrices of rank less than $k$ and $U$ is the subgroup of $GL_k$ equal to block upper triangular matrices with block sizes given by the multiplicities of the $d_i$. In fact, this identification is correct (see \cite{webb}).
\end{example}



\section{A formula for the $I$-function}\label{sec:mainthm}
The quasimap $I$-function was advertised as a more tractable member of the Gromov-Witten family, a generating function that recovers Gromov-Witten invariants after a mirror map, but which (supposedly) has a closed formula. But how, in practice, can one write down the quasimap $I$-function of $V\sslash G$? 

When $G$ is abelian, the answer is well-known (see \cite{givental} or \cite{toric_qmaps}). To handle nonabelian groups, the strategy is to prove an \textit{abelian/nonabelian correspondence} for $I$-functions. ``Abelian/nonabelian correspondences'' can be found throughout representation theory and geometry; the basic idea of these results, if $T$ is a maximal torus of $G$, is to relate some data determined by $G$ to the corresponding data determined by $T$. In our situation, a character $\theta$ of $G$ restricts to a character of $T$, so one may attempt to recover the various invariants (including the $I$-function) of $V\sslash_\theta G$ from the more accessible invariants of $V\sslash_\theta T$. Most abelian/nonabelian correspondences use the \textit{Weyl group} $W$ of $T$ in $G$, which equals the quotient of the normalizer of $T$ by the torus $T$ itself: $W = N_G(T)/T$.

\begin{example}\label{ex:tau}
We describe an ``abelian/nonabelian correspondence'' for principal bundles on $\PP^1$. Let $\T$ be a principal $T$-bundle; then if we let $T$ act on $G$ by left multiplication, the associated bundle $\T\times_T G$ is a principal $G$-bundle. There is a group homomorphism $\tau: \Hom(\chi(T),\ZZ) \rightarrow \Hom(\chi(G),\ZZ)$ induced by restriction of characters, and if $\T$ has degree $\tilde \beta \in \Hom(\chi(T),\ZZ)$, then the degree of $\T\times_T G$ is $\tau(\tilde \beta)$. Moreover, the Weyl group of $T$ in $G$ acts on $\Hom(\chi(T),\ZZ)$. Grothendieck's classification theorem \cite{gro57} says that every principal $G$-bundle may be written $\T\times_T G$ for some $\T$, and that the isomorphism class of $\T\times_T G$ is determined by the Weyl orbit of $\tilde \beta$ in $\Hom(\chi(T),\ZZ)$.
\end{example}

We may understand the cohomology of $V\sslash G$ as follows. There is a rational map $p: V\sslash T \dashrightarrow V\sslash G$ defined on $A=V^s(G)/T \subset V^s(T)/T$. Let $W$ denote the Weyl group of $T$ in $G$; then $W$ acts on $A$ and hence on $H^*(A,\QQ)$. A classical argument shows that $p^*$ defines an isomorphism from $H^*(V\sslash G, \QQ)$ to $H^*(A, \QQ)^W$ (this may be read as an abelian/nonabelian correspondence for cohomology
\footnote{The full relationship of $H^*(V\sslash G, \QQ)$ and $H^*(V\sslash T, \QQ)$ (which requires understanding the restriction $H^*(V\sslash T, \QQ)\rightarrow H^*(A, \QQ)$) is worked out in \cite{martin} by Martin; the analogous result in Chow is proved in \cite{ellingsrud} by Ellingsrud-Stromme.}). This is useful because we can write down explicit classes in $H^*(A,\QQ)$ from characters of $T$: If $\xi \in \chi(T)$, let $L_\xi$ denote the line bundle $V^s(G)\times_{\CC^*}\CC_\xi$ on $A$. This gives a class $c_1(L_\xi)$, and with the right definitions, the map $\xi \mapsto c_1(L_\xi)$ is $W$-equivariant.

With this language, we can now write down a formula for the quasimap $I$-function. Because the formula expresses the $I$-function for $V\sslash G$ in terms of the $I$-function for $V\sslash T$ and the Weyl group, we call it an abelian/nonabelian correspondence for $I$-functions.
\begin{theorem}\cite{webb}\label{thm:mainthm}
Let $G$ be a reductive algebraic group and let $T$ be a maximal torus of $G$. Suppose $G$ acts on a vector space $V$, and that the weights of the restriction of this action to $T$ are $\xi_1, \ldots, \xi_n \in \chi(T)$. Then the coefficients of the $I$-function of $V\sslash G$ satisfy 
\begin{equation}\label{eq:mainresult}
I^{V\sslash G}_{\beta}(z) = \sum_{\tilde \beta \rightarrow \beta} \left( \prod_{\alpha} \frac{\prod_{k=-\infty}^{\tilde \beta \cdot c_1(L_{\alpha})}(c_1(L_{\alpha}) + kz)}{\prod_{k=-\infty}^0 (c_1(L_{\alpha}) + kz)}\right)I^{V\sslash T}_{\tilde \beta}(z) 
\end{equation}
where
\begin{equation}\label{eq:toric}\quad I_{\tilde \beta}^{V\sslash T}(z) = \prod_{i=1}^n \frac{\prod_{k=-\infty}^0(c_1(L_{\xi_i})+kz)}{\prod_{k=-\infty}^{\tilde \beta(\xi_i)}(c_1(L_{\xi_i})+kz)}.
\end{equation}
The sum is over all effective $\tilde \beta$ in $\tau^{-1}(\beta)$, the variable $\alpha$ runs over the roots of $G$, and the equality in \eqref{eq:mainresult} is equality after pulling both sides back to $H^*(A,\QQ)$.
\end{theorem}

The formula \eqref{eq:toric} for the coefficients of a toric $I$-function was originally given by Givental in \cite{givental} (note that the quasimap $I$-function differs from Givental's $I$-function by an exponential factor). Theorem \ref{thm:mainthm} has natural generalizations to twisted and equivariant $I$-functions.  

One application of this theorem is an abelian/nonabelian correspondence in quantum cohomology for a large class of targets. Indeed, for these targets, the correspondence was reduced to a correspondence of small $J$-functions in \cite{frob-ab-nonab}, and the wall-crossing result of \cite{wcgis0} translates this to the relationship of small $I$-functions in Theorem \ref{thm:mainthm}. (A more general abelian/nonabelian correspondence theorem in quantum cohomology was proved using symplectic geometry by Gonzalez-Woodward in \cite{gonzalez}.) A second application, due to Kalashnikov \cite{elana}, is to identify nonisomorphic Fano 4-folds by comparing their quantum periods, which are determined by their $I$-functions (this reference also provides an independent proof of Theorem \ref{thm:mainthm} when the target is a quiver flag variety). Thirdly, the proof of Theorem \ref{thm:mainthm} is very geometric, and so translates easily to the setting of K-theory, providing an analogous formula for the K-theoretic $I$-function. By providing explicit formulas for $I$-functions, Theorem \ref{thm:mainthm} should facilitate the computation of examples in Gromov-Witten theory, and may lead to the discovery of new relationships.

We close this paper with a sketch of the proof of Theorem \ref{thm:mainthm}.

\begin{proof}[Sketch of proof of Theorem \ref{thm:mainthm}] 
To compute the $I$-function coefficients $I^{V\sslash G}_\beta(z)$ in \eqref{eq:Ifunc}, the key step is to identify the fixed locus $F_\beta$. In fact, connected components of $F_\beta$ are indexed by the isomorphism type of the principal bundle $\P$ of the quasimap, or equivalently by the Weyl orbits in $\tau^{-1}(\beta)$ (see Example \ref{ex:tau}). Choose representatives $\tilde \beta_i$ for these orbits and let $F_{\tilde \beta_i}$ represent the corresponding component of the fixed locus. As we did in Example \ref{ex:Fbeta} for Grassmannians, we can na\"ively guess what variety $F_{\tilde \beta_i}$ should be.

Generalizing Example \ref{ex:qmaps}, a quasimap from $\PP^1$ to $V\sslash G$ may be written as a vector of homogeneous polynomials. If it is in $F_\beta$, then it has a unique basepoint at $[0:1],$ so we hope that we can choose the polynomials to be polynomials in $x$ only. Hence, 
if $T$ acts on $V$ with weights $\xi_1, \ldots, \xi_n$, define $V_{\tilde \beta_i}$ to be the vector space with elements
\begin{equation}\label{eq:poly}
(u_1(x), u_2(x), \ldots, u_n(x)),
\end{equation}
where $u_j(x)$ is a constant multiple of $x^{\tilde \beta_i(\xi_j)}$ if $\tilde \beta_i(\xi_j) \geq 0,$ and $u_j(x)$ is zero otherwise. There is an evaluation map $ev_\bullet: V_{\tilde \beta_i} \rightarrow V$ sending \eqref{eq:poly} to $(u_1(1),\ldots,u_n(1))$. This map embeds $V_{\tilde \beta}$ as the subspace of $V$ where $\tilde \beta_i(\xi_j)\geq 0$, which allows us to intersect $V_{\tilde \beta_i}$ with $V^{s}(G)$ and obtain a space $V^{s}_{\tilde \beta_i}$. Our guess is that every closed point of $F_{\tilde \beta_i}$ may be represented by a point of $V^{s}_{\tilde \beta_i}$.

However, distinct points of $V^{s}_{\tilde \beta_i}$ may not represent distinct quasimaps. Generalizing Example \ref{ex:qmap_isos}, an isomorphism of quasimaps in $V^{s}_{\tilde \beta_i}$ is an invertible $n\times n$ matrix of polynomials $(a_{\ell j}(x))_{\ell j}$ where $a_{\ell j}(x)$ is homogeneous of degree $\tilde \beta_i(\xi_\ell)-\tilde \beta_i(\xi_j)$ if this quantity is nonnegative, and $a_{\ell j}(x)=0$ otherwise. Let $P_{\tilde \beta_i}$ denote this group of matrices; we also have an evaluation map $ev_\bullet: P_{\tilde \beta_i} \rightarrow \Aut(V)$ by evaluating a matrix at 1. In fact, the image of $ev_\bullet$ is in $G$, and $ev_\bullet$ identifies $P_{\tilde \beta_i}$ with a parabolic subgroup of $G$.

Finally we notice that $V^{s}_{\tilde \beta_i}/P_{\tilde \beta_i}$ carries a tautological family of quasimaps with the following section:
\begin{equation}
\begin{aligned}\label{eq:universal_family}
U: \frac{V_{\tilde \beta_i}^s \times (\CC^2 \setminus \{0\})}{(\u, \x) \sim (A \u, t\x)} &\longrightarrow \frac{V_{\tilde \beta_i}^s \times (\CC^2\setminus\{0\})\times V}{(\u, \x, \v) \sim (A \u, t\x, [t^{\tilde \beta(\bxi)}]A(\x)\v)}  & \\
(\u, \x) &\mapsto (\u, \x, \u(\x))
\end{aligned}
\end{equation}
where $(\u, \x, \v)$ is in $V^s_{\tilde \beta_i}\times(\CC^2\setminus\{0\})\times V$ and $(A, t) \in P_{\tilde \beta_i} \times \CC^*,$ and $[t^{\tilde \beta_i(\bxi)}]$ denotes the matrix with diagonal $(t^{\tilde \beta_i(\xi_1)}, \ldots, t^{\tilde \beta_i(\xi_n)}).$ This matrix factor guarantees that every geometric fiber of this family of vector bundles has isomorphism type $\tilde \beta_i$. This tautological family defines a map from $V^s_{\tilde \beta_i}/P_{\tilde \beta_i}$ to $F_{\tilde \beta_i}$, and the hardest part of Theorem \ref{thm:mainthm} is showing that this morphism is an isomorphism. The proof uses that obstructions on $F_{\tilde \beta_i}$ vanish, and moreover this stack is smooth.

From here, the computation of the formula in Theorem \ref{thm:mainthm} is relatively straightforward. The evaluation map $V^s_{\tilde \beta_i}/P_{\tilde \beta_i} \rightarrow V\sslash G$ (which sends $x$ to 1) factors as
\[
V^s_{\tilde \beta_i}/P_{\tilde \beta_i} \xrightarrow{i} V^s/P_{\tilde \beta_i} \xrightarrow{h} V^s/G,
\]
where the first map is a closed embedding (recall $V_{\tilde \beta_i}$ is a subspace of $V$) and $h$ is a flag bundle ($P_{\tilde \beta_i}$ is a parabolic subgroup of $G$). Recall that the equality \eqref{eq:mainresult} holds only after pulling back both sides to $V^s(G)/T$, so let $\Psi: V^s(G)/T\rightarrow V^s(G)/G$ be the projection. We compute
\[
\Psi^*I^{V\sslash G}_{\beta}(z) = \Psi^*\sum_{\tilde \beta_i}(ev_\bullet)_*\frac{[F_{\tilde \beta_i}]^{\vir}}{e_{\CC^*}(N_{F_{\tilde \beta_i}}^{\vir} )} = \sum_{\tilde \beta_i}\Psi^*h_*i_*\frac{1}{e_{\CC^*}(N_{F_{\tilde \beta_i}}^{\vir})}
\]
where we have used that $[F_{\tilde \beta_i}]^{\vir}=1$ because $F_{\tilde \beta_i}$ is smooth and the obstructions vanish (see Example \ref{ex:obs_vanish}). 

To compute $e_{\CC^*}(N_{F_{\tilde \beta_i}}^{\vir})$, we use that the restriction of the absolute obstruction theory to $F_{\tilde \beta_i}$ is $R^\bullet \pi_* \mathcal{F}$, where $\mathcal{F}$ is the sheaf on $(V^s_{\tilde \beta_i}/P_{\tilde \beta_i})\times\PP^1$ is defined by the exact sequence
\[
0 \rightarrow \mathcal{A} \rightarrow \mathcal{B} \rightarrow \mathcal{F} \rightarrow 0
\]
where
{\small
\[
\mathcal{A} = \frac{V^s_{\tilde \beta_i} \times (\CC^2\setminus\{0\})\times \g}{(\u, \x, X) \sim (A\u, t\x, [t^{\tilde \beta_i(\bxi)}]A(\x) \cdot X)}
\quad \text{and}\quad
\mathcal{B} = \frac{V^s_{\tilde \beta_i} \times (\CC^2\setminus\{0\})\times V}{(\u, \x, v) \sim (A\u, t\x, [t^{\tilde \beta_i(\bxi)}]A(\x)\cdot v)}.
\]
}
So $e_{\CC^*}(N_{F_{\tilde \beta_i}}^{\vir})$ is a certain ratio of the moving parts of the pushforwards of $\mathcal{A}$ and $\mathcal{B}$ (see \eqref{eq:nextstep}). We compute $i_*$ with the projection formula, noting that $V^s_{\tilde \beta_i}$ is the zero locus of the tautological section of the bundle $V^s \times (V/V_{\tilde \beta_i})$ on $V^s$. We compute $\Psi^*h_*$ with a lemma of Brion \cite{brion96}, which introduces an additional sum over the Weyl group. After these steps we have
\begin{equation}\label{eq:nextstep}
\Psi^*I^{V\sslash G}_{\beta}(z) = \sum_{\tilde \beta_i} \sum_{\omega \in W/W_{L_i}} \omega\left[\frac{e(V^s \times_T (V/V_{\tilde \beta_i})) \;\;e_{\CC^*}(R^0\pi_*\mathcal{A})^{\mov}e_{\CC^*}(R^1\pi_*\mathcal{B})^{\mov}}{\prod_{\alpha \in R^+\setminus R_{L_i}} c_1(L_\alpha)\;\;e_{\CC^*}(R^0\pi_*\mathcal{B})^{\mov}e_{\CC^*}({R^1\pi_*\mathcal{A}})^{\mov}}\right]
\end{equation}
where $\mathcal{A}$ and $\mathcal{B}$ now denote the analogous bundles on $V^s(G)/T$, and $L_i$ is the Levi subgroup of $P_{\tilde \beta_i}$ containting $T$ with roots $R_{L_i}$ and Weyl group $W_{L_i}$, and $R^+$ are the opposite roots of a Borel subgroup of $G$ contained in $P_{\tilde \beta_i}$. 

One can check that $W_{L_i}$ is precisely the stabilizer of $\tilde \beta_i$ in $\Hom(\chi(T),\ZZ)$ and that the double sum and group action may be simplified to a single sum over all $\tilde \beta$ in $\tau^{-1}(\beta)$, giving
\[
\Psi^*I^{V\sslash G}_{\beta}(z) = \sum_{\tilde \beta \rightarrow \beta} \frac{e(V^s \times_T (V/V_{\tilde \beta})) \;\;e_{\CC^*}(R^0\pi_*\mathcal{A})^{\mov}e_{\CC^*}(R^1\pi_*\mathcal{B})^{\mov}}{\prod_{\alpha \in R^+\setminus R_{L_{\tilde \beta}}} c_1(L_\alpha)\; \;e_{\CC^*}(R^0\pi_*\mathcal{B})^{\mov}e_{\CC^*}({R^1\pi_*\mathcal{A}})^{\mov}}.
\]
From here it remains to compute the euler classes. One finds that the classes coming from $\mathcal{A}$ and the flag pushforward $h_*$ yield the twist factor in \eqref{eq:mainresult}:
\[
\frac{e_{\CC^*}(R^0\pi_*\mathcal{A})^{\mov}}{\prod_{\alpha \in R^+\setminus R_{L_{\tilde \beta}}} c_1(L_\alpha)\;e_{\CC^*}({R^1\pi_*\mathcal{A}})^{\mov}} = \left( \prod_{\alpha} \frac{\prod_{k=-\infty}^{\tilde \beta \cdot c_1(L_{\alpha})}(c_1(L_{\alpha}) + kz)}{\prod_{k=-\infty}^0 (c_1(L_{\alpha}) + kz)}\right)
\]
while the classes coming from $\mathcal{B}$ and the inclusion pushforward $i_*$ yield the toric $I$-function:
\[
\frac{e(V^s \times_T (V/V_{\tilde \beta}))\;e_{\CC^*}(R^1\pi_*\mathcal{B})^{\mov}}{e_{\CC^*}(R^0\pi_*\mathcal{B})^{\mov}} = I_{\tilde \beta}^{V\sslash T}(z).
\]

\end{proof}

\bibliographystyle{amsplain}
\bibliography{references}

\end{document}